\documentclass[11pt]{article}

\usepackage{geometry,graphicx,latexsym,amssymb,verbatim, multicol}
\usepackage{tikz}
\usetikzlibrary{calc}

\newtheorem{thm}{Theorem}
\newtheorem{lem}[thm]{Lemma}

\newtheorem{prob}{Problem}
\newtheorem{cor}[thm]{Corollary}
\newtheorem{ob}[thm]{Observation}
\newtheorem{prop}[thm]{Proposition}

\newcommand{\dontshow}[1]{}

\newcommand{\spn}{{\rm span}}
\newcommand{\rep}{{\rm rep}}

\newcommand{\fd}{{\rm fd}}

\newcommand{\kfd}{{\rm fd}_k}
\newcommand{\ofd}{{\rm fd}_1}

\newcommand{\outr}{\xi_{\rm or}}

\newcommand{\proof}{\noindent\textbf{Proof. }}
\newcommand{\smallqed}{{\tiny ($\Box$)}}
\newcommand{\qed}{$\Box$}

\newcommand{\coro}{{\rm cor}}

\newcommand{\mod}{{\rm mod}}

\newcommand{\barK}{{\overline{K}}}
\newcommand{\bard}{{\overline{d}}}

\def \nH {n_{H}}
\def \nF {n_{F}}

\def \mH {m_{H}}

\def \mG {m_{G}}

\newcommand{\barG}{\overline{G}}

\newcommand{\2}{ \vspace{0.2cm} }

\newcommand{\cG}{{\cal G}}

\newenvironment{unnumbered}[1]{\trivlist \item [\hskip \labelsep {\bf
#1}]\ignorespaces\it}{\endtrivlist}

\def\vertex(#1){\put(#1){\circle*{2}}}
\def\vertexo(#1){\put(#1){\circle{2}}}
\def\vert(#1){\put(#1){\circle*{1.5}}}
\def\verto(#1){\put(#1){\circle{1.5}}}
\def\lab(#1)#2{\put(#1){\makebox(0,0)[c]{#2}}}
\begin{document}




\begin{center}
{\LARGE Fair Domination in Graphs}
\mbox{}\\[8ex]

\begin{multicols}{2}

Yair Caro\\[1ex]
{\small Dept. of Mathematics and Physics\\
University of Haifa-Oranim\\
Tivon 36006, Israel\\
yacaro@kvgeva.org.il}

\columnbreak

Adriana Hansberg\\[1ex]
{\small Dept. de Matem\`atica Aplicada III\\
UPC Barcelona\\
08034 Barcelona, Spain\\
adriana.hansberg@upc.edu}\\[2ex]

\end{multicols}

Michael Henning\\[1ex]
{\small Dept. of Mathematics\\
University of Johannesburg\\
Auckland Park 2006, South Africa\\
mahenning@uj.ac.za}

\mbox{}\\[3ex]

\end{center}

\begin{abstract}
A fair dominating set in a graph $G$ (or FD-set) is a dominating set
$S$ such that all vertices not in $S$ are dominated by the same
number of vertices from $S$; that is, every two vertices not in $S$
have the same number of neighbors in $S$. The fair domination number,
$\fd(G)$, of $G$ is the minimum cardinality of a FD-set. We present
various results on the fair domination number of a graph. In
particular, we show that if $G$ is a connected graph of order~$n \ge
3$ with no isolated vertex, then $\fd(G) \le n - 2$, and we construct
an infinite family of connected graphs achieving equality in this
bound. We show that if $G$ is a maximal outerplanar graph, then
$\fd(G) < 17n/19$. If $T$ is a tree of order~$n \ge 2$, then we prove
that $\fd(T) \le n/2$ with equality if and only if $T$ is the corona
of a tree.
\end{abstract}


\section{Introduction}

In this paper, we continue the study of domination in graphs.
Domination in graphs is now well studied in Graph Theory. The
literature on this subject has been surveyed and detailed in the two
books by Haynes, Hedetniemi, and Slater~\cite{hhs1, hhs2}.
For notation and Graph Theory terminology we in general
follow~\cite{hhs1}. Specifically, let $G = (V, E)$ be a  graph with
vertex set $V$ of order~$n = |V|$ and edge set $E$ of size~$m =
|E|$, and let $v$ be a vertex in $V$. The \emph{open neighborhood}
of $v$ is the set $N(v) = \{u \in V \, | \, uv \in E\}$, while the
\emph{closed neighborhood} of $v$ is the set $N[v] = N(v) \cup
\{v\}$.

Let $G = (V,E)$ be a graph. A \emph{dominating set} in $G$ is a set
$D$ of vertices of $G$ such that every vertex $v \in V$ is either in
$D$ or adjacent to a vertex of $D$. A vertex in $D$ is said to
\emph{dominate} a vertex outside $D$ if they are adjacent in $G$. The
\emph{domination number} of $G$, denoted $\gamma(G)$, is the minimum
cardinality of a dominating set. A dominating set of $G$ of
cardinality $\gamma(G)$ is called a $\gamma(G)$-set.

Let $G$ be a graph that is not the empty graph. For $k \ge 1$ an
integer, a \emph{$k$-fair dominating set}, abbreviated kFD-set, in
$G$ is a dominating set $D$ such that $|N(v) \cap D| = k$ for every
vertex $v \in V \setminus D$. We note that the set $D = V$ is a
kFD-set since vacuously every vertex in $V \setminus D = \emptyset$
satisfies the desired property. The \emph{$k$-fair domination number}
of $G$, denoted by $\kfd(G)$, is the minimum cardinality of a
kFD-set. A kFD-set of $G$ of cardinality $\kfd(G)$ is called a
$\kfd(G)$-set. With this definition in mind, we point to a related
problem on the so called $(k, \tau)$-regular sets discussed in
\cite{CKV} (and all the references given there). Some reminiscent of
this approach appears more explicitly in Proposition \ref{l:line1}.

A \emph{fair dominating set}, abbreviated FD-set, in $G$ is a kFD-set
for some integer $k \ge 1$. Thus a dominating set $D$ is a FD-set in
$G$ if $D = V$ or if $D \ne V$ and all vertices not in $D$ are
dominated by the same number of vertices from $D$; that is, $|N(u)
\cap D| = |N(v) \cap D| > 0$ for every two vertices $u, v \in V
\setminus D$.
We remark that if $G \ne \barK_n$, then $G$ contains a vertex $v$
that is not isolated in $G$ and the set $V \setminus \{v\}$ is a
FD-set in $G$. Hence every graph that is not empty has a FD-set of
cardinality strictly less than its order.
The \emph{fair domination number}, denoted by $\fd(G)$, of a graph
$G$ that is not the empty graph is the minimum cardinality of a
FD-set in $G$. By convention, if $G = \barK_n$, we define $\fd(G) =
n$. Hence if $G$ is not the empty graph, then $\fd(G) = \min \{
\kfd(G) \}$, where the minimum is taken over all integers $k$ where
$1 \le k \le |V|-1$. A FD-set of $G$ of cardinality $\fd(G)$ is
called a $\fd(G)$-set. Every FD-set in a graph $G$ is a dominating
set in $G$. Hence we have the following observation.

\begin{ob}
Let $G$ be a graph of order~$n$. Then the following holds. \\
\indent {\rm (a)} $\gamma(G) \le \fd(G)$. \\
\indent {\rm (b)} $\fd(G) \le n$, with equality if and only if $G =
\barK_n$. \label{ob1}
\end{ob}

We show later (see Corollary~\ref{c:bound}) that the result in
Observation~\ref{ob1}(b) can be improved as follows: if $G$ is a
graph of order~$n$, then $\fd(G) \le n-2$, unless $G = \barK_n$, in
which case $\fd(G) = n$, or $G$ contains precisely one edge, in which
case $\fd(G) = n-1$.

\begin{figure}[htb]
\tikzstyle{every node}=[circle, draw, fill=black!0, inner sep=0pt,minimum width=.2cm]
\begin{center}
\begin{tikzpicture}[thick,scale=.7]
  \draw(0,0) { 
    +(1.90,3.62) -- +(0.00,2.24)
    +(0.00,2.24) -- +(0.73,0.00)
    +(0.73,0.00) -- +(3.08,0.00)
    +(3.08,0.00) -- +(3.80,2.24)
    +(1.90,3.62) -- +(3.80,2.24)
    +(1.90,3.62) -- +(1.90,2.62)
    +(0.00,2.24) -- +(0.95,1.93)
    +(1.31,0.81) -- +(0.73,0.00)
    +(2.49,0.81) -- +(3.08,0.00)
    +(2.85,1.93) -- +(3.80,2.24)
    +(1.90,2.62) -- +(2.49,0.81)
    +(2.49,0.81) -- +(0.95,1.93)
    +(0.95,1.93) -- +(2.85,1.93)
    +(2.85,1.93) -- +(1.31,0.81)
    +(1.31,0.81) -- +(1.90,2.62)
    +(0.95,1.93) node{}
    +(0.00,2.24) node{}
    +(1.90,2.62) node{}
    +(1.90,3.62) node{}
    +(2.85,1.93) node{}
    +(2.49,0.81) node{}
    +(1.31,0.81) node{}
    +(0.73,0.00) node{}
    +(3.08,0.00) node{}
    +(3.80,2.24) node{}
  };
\end{tikzpicture}
\end{center}
\vskip -0.6 cm \caption{The Petersen graph $G_{10}$ with $\fd(G_{10}) = 4$.} \label{Peter}
\end{figure}
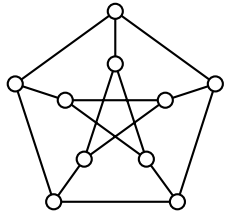

For example, consider the Petersen graph $G = G_{10}$ shown in
Figure~\ref{Peter} which has domination number $\gamma(G) = 3$. Let
$V = V(G)$. The only possible $\gamma(G)$-sets are the open
neighborhoods $D = N(v)$, where $v \in V$, but these are not FD-sets
since if $u \in V \setminus D$, then either $u \ne v$, in which case
$|N(u) \cap D| = 1$, or $u = v$, in which case $|N(u) \cap D| = 3$.
Thus, by Observation~\ref{ob1}(a), $\fd(G) > \gamma(G) = 3$. However
the closed neighborhood $N[v]$ of any vertex $v \in V$ forms a FD-set
in $G$, and so $\fd(G) \le |N[v]| = 4$. Consequently, $\fd(G) = 4$.

An \emph{out-regular set}, abbreviated OR-set, of a graph $G$ that is
not the empty graph is a set $Q$ of vertices such that $|N(u) \cap (
V \setminus Q )| = |N(v) \cap ( V \setminus Q)| > 0$ for every two
vertices $u$ and $v$ in $Q$. We remark that if $G \ne \barK_n$, then
$G$ contains a vertex $v$ that is not isolated in $G$ and the set
$\{v\}$ is an OR-set in $G$. Hence every graph that is not empty has
an OR-set. The \emph{out-regular number} of a non-empty graph $G$,
denoted by $\outr(G)$, is the maximum cardinality of an OR-set. By
convention, if $G = \barK_n$, we define $\outr(G) = 0$. An OR-set of
$G$ of cardinality $\outr(G)$ is called an $\outr(G)$-set. Since
every graph that is not empty has an OR-set, we have the following
observation.

\begin{ob}
Let $G$ be a graph of order~$n$. Then, $\outr(G) \ge 0$, with
equality if and only if $G = \barK_n$. \label{ob2}
\end{ob}

If $S$ is a packing in $G$ (and so, the vertices in $S$ are pairwise
at distance at least~$3$ apart in $G$), then in order to dominate the
vertices in $S$ every dominating set in $G$ must contain at least one
vertex in $N[v]$ for each $v \in S$, and so $\gamma(G) \ge |S|$. A
\emph{perfect dominating set}, abbreviated PD-set, in $G$ is a
dominating set $S$ that is a packing in $G$. Thus if $S$ is a PD-set
in $G$, then every vertex not in $S$ is dominated by a unique vertex
in $S$, and so $S$ is a 1FD-set, implying that $\fd(G) \le \ofd(G)
\le \gamma(G)$. Consequently, by Observation~\ref{ob1}, we have the
following observation.

\begin{ob}
If a graph $G$ has a PD-set, then $\gamma(G) = \ofd(G)  = \fd(G)$. \label{ob3}
\end{ob}

\subsection{Notation}

%
%
%
We denote the \emph{degree} of $v$ in $G$ by $d_G(v)$, or simply by
$d(v)$ if the graph $G$ is clear from context. Let $\delta(G)$,
$\Delta(G)$ and $\bard(G)$ denote, respectively, the minimum degree,
the maximum degree and the average degree in $G$. For a set
$S\subseteq V$, we denote the number of vertices of $S$ adjacent to
$v$ in $G$ by $d_S(v)$. In particular, $d_V(v) = d_G(v)$. For a set
$S \subseteq V$, the subgraph induced by $S$ is denoted by $G[S]$. We
denote by $\spn(G)$, the \emph{span} of $G$, the number of distinct
values in the degree sequence of $G$ and by $\rep(G)$, the
\emph{repetition number} of $G$, the maximum multiplicity in the list
of vertex degrees.
The parameter $\alpha(G)$ denotes the (vertex) \emph{independence
number} of $G$, while $\chi(G)$ denotes the \emph{chromatic number}
of $G$.

Further, we denote the complete graph on $n$ vertices by $K_n$ and
the empty graph on $n$ vertices by $\barK_n$. Moreover, $P_n$, $C_n$
and $K_{m,n}$ denote, respectively, the path on $n$ vertices, the
cycle on $n$ vertices and the complete bipartite graph with one
partite set of cardinality~$m$ and the other of cardinality~$n$.


\section{Preliminary Results and Observations}

It is a simple exercise to determine the fair domination number of
certain well-studied families of graphs. Recall that for $n \ge 3$,
$\gamma(P_n) = \gamma(C_n) = \lceil n/3 \rceil$, while for $m \ge n
\ge 2$, $\gamma(K_{m,n}) = 2$.

\begin{ob} \label{ob4}
For $m, n \ge 1$, if $G \in \{P_n,K_n,\barK_n,K_{m,n}\}$, then
$\fd(G) = \gamma(G)$. Further for $n \ge 3$, $\fd(C_n) = \gamma(C_n)$
unless $n \equiv 2 \, (\mod \, 3)$ and $n \ge 5$ in which case
$\fd(C_n) = \gamma(C_n) + 1$.
\end{ob}

We next establish a relationship between the fair domination number
and the out-regular number of a graph.

\begin{prop}
For every graph $G$ of order~$n \ge 2$, $\fd(G) + \outr(G) = n$.
\label{prop1}
\end{prop}
\textbf{Proof.} If $G = \barK_n$, then $\fd(G) = n$ and, by
convention, $\outr(G) = 0$. Hence we may assume that $G \ne \barK_n$,
for otherwise the desired result holds. Let $D$ be a $\fd(G)$-set. By
Observation~\ref{ob1}(b), $\fd(G) < n$. Let $Q = V \setminus D$.
Then, $Q$ is an OR-set in $G$, and so $\outr(G) \ge |Q| = n -
\fd(G)$, or, equivalently, $\fd(G) + \outr(G) \ge n$. Conversely, let
$Q$ be an $\outr(G)$-set. By Observation~\ref{ob2}, $\outr(G)
> 0$. By definition, $\outr(G) < n$. Let $D = V \setminus Q$. Then,
$D$ is a FD-set, and so $\fd(G) \le |D| = n - \outr(G)$, or,
equivalently,  $\fd(G) + \outr(G) \le n$. Consequently, $\fd(G) +
\outr(G) = n$.~\qed

\begin{thm}
Let G be a connected graph on $n \ge 2$ vertices. Then the following
holds. \\
\indent {\rm (a)} If $\barG$ is connected, then $\fd(G) =
\fd(\barG)$. \\
\indent {\rm (b)} If $\barG$ has~$q \ge 2$ components, then $\fd(G)
\le n/q \le n/2$. \label{thm:Gbar}
\end{thm}
\textbf{Proof.} (a) Suppose that $\barG$ is connected. Let $D$ be a
$\fd(G)$-set. Then every vertex $v \in V \setminus D$ is adjacent to
precisely~$k$ vertices in $D$ for some integer $k$, $1 \le k \le
|D|$. If $k = |D|$, then in $\barG$ there are no edges between $D$
and $V \setminus D$, contradicting the assumption that $\barG$ is
connected. Hence, $k < |D|$. But then in $\barG$ every vertex in $V
\setminus D$ is adjacent to precisely $|D| - k > 0$ vertices in $D$,
and so $D$ is a FD-set in $\barG$. Thus, $\fd(\barG) \le |D| =
\fd(G)$. Reversing the roles of $G$ and $\barG$, we have that $\fd(G)
\le \fd(\barG)$. Consequently, $\fd(G) = \fd(\barG)$.

(b) Suppose that $\barG$ is not connected and has~$q$ components. Clearly, the smallest component in $\barG$ has
cardinality at most~$n/q$. Let $F$ be the smallest component in
$\barG$ and let $D = V(F)$. Then in $G$ every vertex in $V \setminus
D$ is adjacent to all vertices in $D$, and so $D$ is a FD-set in $G$.
Thus, $\fd(G) \le |D| \le n/q \le n/2$.~\qed

\medskip
Next we consider the fair domination number of the line graph,
$L(G)$, of a graph $G$.

\begin{prop}
Let $G$ be a graph of size~$\mG$ and let $L(G)$ denote the line graph
of $G$. If $H$ is a spanning $r$-regular subgraph of $G$, where $r
> 0$ and where $H$ is not necessarily induced, of size~$\mH$, then
\[
\fd(L(G)) \le \mH = \left( \frac{\mH}{\mG} \right) |V(L(G)))|.
\]
\label{l:line1}
\end{prop}
\textbf{Proof.} Let $D$ be the set of vertices in the line graph
$L(G)$ of $G$ corresponding to the edges in $H$, and let $v \in
V(L(G)) \setminus D$. Thus the vertex $v$ corresponds to an edge $e
\in E(G) \setminus E(H)$, and so both ends of $e$ are incident with
precisely~$r$ edges of $H$. Hence in $L(G)$, the vertex $v$ is
adjacent to exactly~$2r$ vertices in $D$. Thus in $L(G)$, we have
that $|N(v) \cap D| = 2r$ for every vertex $v \in V(L(G)) \setminus
D$, implying that $D$ is a FD-set in $L(G)$, and so $\fd(L(G)) \le
|D| = \mH$.~\qed

\medskip
We remark that examples of graphs that possess spanning $r$-regular
subgraphs, where $r > 0$ and where $H$ is not necessarily induced,
are abundant. For example, regular graphs of even degree have a
$2$-factor as do Hamiltonian graphs. Several interesting families of
graphs possess a $1$-factor (or perfect matching), including regular
bipartite graphs and connected claw-free graphs of even order. For
further results about regular spanning graphs see for example
\cite{Cs07,HoVo}.

\section{Results}

\subsection{Upper Bounds}

We first establish upper bounds on the fair domination of a graph in
terms of its order. By Observation~\ref{ob3}(a), $\fd(G) \le n$ with
equality if $G = \barK_n$. However this bound can be improved
slightly if we restrict our attention to graphs without isolated
vertices.

\begin{thm}
If $G$ is a graph of order~$n \ge 3$ with $\delta(G) \ge 1$, then
$\fd(G) \le n - 2$, and this bound is sharp. \label{thm_n2}
\end{thm}
\textbf{Proof.} We proceed by induction on $n \ge 3$. If $n \in
\{3,4\}$, it is a simple case to check that if $G$ is a graph of
order~$n$ with no isolated vertex, then $\fd(G) \le n - 2$. This
establishes the base cases. Let $n \ge 5$ and assume that every graph
$G'$ of order~$n'$, where $3 \le n' < n$, with no isolated vertex
satisfies $\fd(G') \le n' - 2$. Let $G = (V,E)$ be a graph of
order~$n$ with no isolated vertex. Every graph on at least two
vertices has two vertices of the same degree. Let $u$ and $v$ be two
vertices in $G$ with the same degree. Then, $d_G(u) = d_G(v) = k$ for
some $k$, where $1 \le k \le n-1$. If $u$ and $v$ are not adjacent in
$G$ or if $u$ and $v$ are adjacent in $G$ and $k \ge 2$, then $D = V
\setminus \{u,v\}$ is a FD-set, and so $\fd(G) \le |D| = n - 2$, as
desired. Hence we may assume that $u$ and $v$ are adjacent in $G$ and
$k = 1$. We now consider the graph $G' = G - \{u,v\}$ of order~$n' =
n - 2 \ge 3$. Since $G$ has no isolated vertex, neither does $G'$.
Applying the inductive hypothesis to $G'$, we have that $\fd(G') \le
n' - 2$. Every $\fd(G')$-set can be extended to a FD-set in $G$ by
adding to it the vertices $u$ and $v$, implying that $\fd(G) \le
\fd(G') + 2 \le n' = n-2$. This establishes the desired upper bound.

We show next that the upper bound is sharp. For this purpose we
construct an infinite family of graphs $G$ of order~$n \ge 3$ with
$\delta(G) \ge 1$ satisfying $\fd(G) = n - 2$. We consider two cases
in turn, depending on the parity of $n$.

\begin{unnumbered}{Claim~I.} There exists an infinite family of
graphs $H$ of even order~$\nH$ satisfying $\fd(H) = \nH - 2$.
\end{unnumbered}
\textbf{Proof.} For $n \ge 3$ define the graph $H = H_n$ on $\nH =
2n$ vertices as follows: Let $V(H) = X \cup Y$, where $X =
\{x_1,\ldots,x_n\}$ and $Y = \{y_1,\ldots,y_n\}$, and where $x_i$ is
adjacent to $y_j$ if and only if $i \ge j$. Further $Y$ is an independent set and, for $i,j > 1$,
$x_i$ is adjacent to $x_j$. Thus, $d_H(x_1)  = 1$ and $d_H(x_i) = i +
n - 2$  for $1 < i \le n$, while $d_H(y_i)  = n- i +1$ for $1 \le i
\le n$. We note that $d_H(x_1) = d_H(y_n) = 1$ and $d_H(y_1) =
d_H(x_2) = n$ and the degrees of all other vertices in $H$ are
distinct.
We show that $\fd(H) = 2n-2 = \nH - 2$. Let $S$ be an arbitrary
FD-set of $H$. We show that $|S| \ge \nH-2$. We consider three cases
in turn.

\indent \emph{Case~I.1. $x_1 \notin S$}. Then evidently $y_1 \in S$
and $S$ is a 1FD-set of $H_n$. If $x_i \notin S$ for some $i \ge 2$,
then $N(x_i) \cap S = \{y_1\}$ and thus $y_2 \notin S$. Since $y_2$
has to be dominated by $S$, there is some $l \ge 2$, $l \neq i$, such
that $x_l \in S$, implying that $y_1, x_l \in N(x_i) \cap S$, a
contradiction. Hence, $X \setminus \{x_1\} \subseteq S$. Moreover,
$y_j \in S$ for any $1 < j < n$, since $|N(y_j) \cap S| \ge |N(y_j)
\cap (X \setminus \{x_1\})| = n-j+1 \ge 2$ and $S$ is a 1FD-set. Thus
$V \setminus \{x_1, y_n\} \subseteq S$ and $|S| \ge \nH-2$, so we are
done.

\indent \emph{Case~I.2. $x_1 \in S$ and $y_n \notin S$.}  Then
evidently $x_n \in S$ and $S$ is a 1FD-set of $H_n$. Since $x_1, x_n
\in N(y_1) \cap S$, it follows that $y_1 \in S$. Hence, as $x_n, y_1
\in N(x_i) \cap S$  for any $1< i < n$, we have that $X \subseteq S$.
This implies also that $y_j \in S$ for any $1 <j < n$ and thus $V
\setminus \{y_n\} \subseteq S$ and therefore $|S| \ge \nH-1$.

\indent \emph{Case~I.3. $x_1, y_n \in S$.} We divide this case in
three parts.

\noindent \emph{(i) Suppose that $y_1 \notin S$.} Then $y_j \in S$
for every $1 < j < n$, otherwise $|N(y_1) \cap S| > |N(y_j) \cap S|$
since $N(y_j) \subset N(y_1)$ and $x_1 \in N(y_1) \setminus N(y_j)$.
As $N(x_i) \varsubsetneqq N(x_{i+1})$ for all $1 < i < n-1$, it
follows that there is at most one index $1 < i < n$ such that $x_i
\notin S$, implying that $|S| \ge \nH-2$.

\noindent \emph{(ii) Suppose that $y_1 \in S$ and $x_n \notin S$.}
Then as $N(x_i) \subset N(x_n)$ and $y_n \in S \cap (N(x_n) \setminus
N(x_i))$ for each $1 < i < n$, it follows that $X \setminus \{x_n\}
\subset S$. Hence, as every $y_j$ has a different number of neighbors
in $X \setminus \{x_n\}$, $y_j \in S$ for all $1 < j < n$. Hence $|S|
\ge \nH-1$.

\noindent \emph{(iii) Suppose that $y_1 \in S$ and $x_m \in S$.} If
$x_i \notin S$ for some $1 < i < m$, then $y_j \in S$ for every $1 <
j < m$, since $N(y_j) \subset N(x_i)$ and $y_1 \in S \cap (N(x_i)
\setminus N(y_j))$. Then $X \setminus \{x_i\} \subset S$ and thus
$|S| \ge \nH-1$ and we are done. Thus assume that $X \subset S$. Then
$|Y \cap S| \ge n-1$, otherwise two different vertices from $Y$ would
not be fairly dominated by $S$. Hence again $|S| \ge \nH-1$.

In all three cases, we have that $|S| \ge \nH-2$. Since $S$ is an
arbitrary FD-set of $H$, it follows that $\fd(H) \ge \nH - 2$. However
as shown earlier, $\fd(H) \le \nH - 2$. Consequently, $\fd(H) = \nH -
2$. This completes the proof of Claim~I.~\smallqed

\begin{unnumbered}{Claim~II.} There exists an infinite family of
graphs $F$ of odd order~$\nF$ satisfying $\fd(F) = \nF - 2$.
\end{unnumbered}
\textbf{Proof.} For $n \ge 3$ define the graph $F = F_n$ on $\nF =
2n+1$ vertices as follows: Let $F$ be obtained from the graph $H_n$
of even order $\nH = 2n$ defined in Claim~I by adding a new vertex
$x_{n+1}$ and joining it to every vertex in $X \setminus \{x_1\}$.
Let $X_F = X \cup \{x_{n+1}\}$. For the proof, consider an arbitrary
FD-set of $F$ and show that $|S| \ge \nF-2$. The proof to show that
for an arbitrary FD-set $S$ of $F$ we have $|S| \ge \nF-2$ is similar
to the proof presented in Claim~I and is therefore omitted.~\smallqed

\medskip
By Claim~I and Claim~II, there is an infinite family of graphs $G$ of
order~$n \ge 6$ with $\delta(G) \ge 1$ satisfying $\fd(G) = n - 2$,
irrespective of whether $n$ is even or odd.~\qed

\medskip
Note that, by Claims~I and II, for each integer $n \ge 6$, there is a
connected graph $G$ on $n$ vertices satisfying $\fd(G) = n - 2$.  If
$n \in \{3,4,5\}$, we can simply take $G = C_n$. Hence for all $n \ge
3$, there exists a connected graph $G$ on $n$ vertices satisfying
$\fd(G) = n - 2$.

\medskip
Let $G$ be a graph of order~$n$ with at least two edges. Let $G^*$ be
the subgraph obtained from $G$ by deleting all isolated vertices in
$G$, if any. Then, $\delta(G^*) \ge 1$ and $G^*$ has order~$n^* \ge
3$. Applying Theorem~\ref{thm_n2} to $G^*$, we have that $\fd(G^*)
\le n^*-2$. Since every $\fd(G^*)$-set can be extended to a FD-set in
$G$ by adding to it the set of isolated vertices in $G$, we have that
$\fd(G) \le n-2$. Hence as a consequence of Theorem~\ref{thm_n2}, we
have the following result.

\begin{cor}
If $G$ is a graph of order~$n$ and size at least~$2$, then $\fd(G)
\le n-2$. \label{c:bound}
\end{cor}

We next present some upper bounds on the fair domination number in terms
of its order, chromatic number and average, maximum and minimum degrees. For this purpose,
we first recall the Caro-Wei Theorem (see~\cite{Ca79,Wei81}).

\begin{unnumbered}{Caro-Wei Theorem.} For every graph $G$ of
order~$n$, \[
\alpha(G) \ge \sum_{v \in V(G)} \frac{1}{1 + d_G(v)} \ge
\frac{n}{\bard(G) + 1}.
\]
\end{unnumbered}

\medskip
For our purposes, we also need the following useful lower bound on
the repetition number of a graph, established by Caro and West in
\cite{CaWe09}.

\begin{unnumbered}{Caro-West Lemma.}
If $G$ is a graph of order~$n$, then $\rep(G) \ge n/(2\bard(G)
-2\delta(G) + 1)$.
\end{unnumbered}

\begin{prop}
Let $G$ be a graph of order~$n$. Then the following holds. \\
{\rm (a)} If $n \ge 2$ and $\delta(G) \ge 1$, then $\fd(G) \le n - n/((\bard(G) + 1) \Delta(G))$.\\
{\rm (b)} $\fd(G) \le n  - n/(2\bard(G) - 2\delta(G) +1)\chi(G)$. \\
{\rm (c)} For $r \ge 2$, if $G$ is an $r$-regular graph, then $\fd(G)
\le rn/(r+1)$. \label{prop2}
\end{prop}
\textbf{Proof.} (a) Let $B$ be a maximum independent set in $G$. Since
$G$ has no isolated vertex, the set $V \setminus B$ is a dominating
set in $G$. We now consider the degrees of the vertices of $B$. There
are at most $\spn(G)$ possible distinct values for these degrees, and
so, by the Pigeonhole Principle, at least one value, say $q$, appears
at least $\alpha(G)/\spn(G)$ times. Let $Q$ be the set of all
vertices in $B$ with degree~$q$. Then, $Q \subseteq B$ and $|Q| = q$.
Let $D = V \setminus Q$. Then, $D$ is a FD-set, and so $\fd(G) \le n
- |Q| \le n - \alpha(G)/\spn(G)$. The desired result now follows from
the observation that $\spn(G) \le \Delta(G)$ and from the Caro-Wei
Theorem.

 (b)  Let $\rep(G) = m$ and suppose that $X =
\{v_1,\ldots,v_m\}$ is a set of vertices with the same degree in $G$.
Let $H = G[X]$ be the subgraph induced by the set $X$.  Clearly,
$\chi(H) \le \chi(G)$, and so $\alpha(H) \ge m/\chi(H) \ge m/\chi(G)
\ge n/(2\bard(G) -2\delta(G) +1)\chi(G)$ by the Caro-West Lemma. Let
$Q$ be a maximum independent set in $H$, and so $|Q| = \alpha(H)$.
Let $D = V(G) \setminus Q$. Then, $D$ is a FD-set in $G$, and so
$\fd(G) \le n - |Q| \le n - \alpha(H)$, and the desired result
follows.

(c) Since $G$ is an $r$-regular graph, we note that $\bard(G) =
\delta(G)$ and $\chi(G) \le r+1$ and the result follows from
Part~(b).~\qed

\medskip
Note that, in case that $G$ is a regular graph, the complement of a
fair dominating set of $G$ is an induced regular subgraph of $G$. At
this point, it is worth mentioning the famous
Erd\H{o}s-Fajtlowicz-Staton problem (see \cite{AlKrSu}) about the
largest induced regular subgraph of a graph $G$. In the case when $G$
is regular, the complement of such a subgraph is a fair dominating
set of $G$ .

\begin{prop}
For $r \ge 1$, if $G$ is an $r$-regular graph on $n$ vertices, then
$\fd(G) \le n - c \log n$ for some $c > 0$. \label{prop4}
\end{prop}
\textbf{Proof.} Let $\cG_n$ denote the family of all graphs of
order~$n$. By Ramsey's theory, for all graphs $G \in \cG_n$ we have
$\max \{ \alpha(G), \alpha(\barG) \} \ge c \log n$ for some constant
$c$. For $r \ge 1$, let $G$ be an $r$-regular graph in $\cG_n$. Then,
$G$ contains either an independent set or a clique of order at
least~$c\log n$ for some constant $c$. Let $X$ be the vertex set of
such an independent set or clique in $G$. Then, $|X| \ge c\log n$ and
the subgraph $G[X]$ induced by $X$ is $s$-regular for some $s$. If
$X$ is an independent set in $G$, then $s = 0$, while if $X$ is a
clique, then $s = |X| - 1 \le r$. Further if $|X| - 1 = r$, then we
note that $K_{r+1}$ is a component of $G$. On the one hand, if $s <
r$, then for every vertex $v \in X$,  we have $|N(v) \cap (V
\setminus X)| = r-s > 0$ and hence $V \setminus X$ is a FD-set of
$G$, and so $\fd(G) \le n - |X| \le n - c\log n$. On the other hand,
if $s = r$, then we choose a vertex $v \in X$ and define $D = (V
\setminus X) \cup \{v\}$. Then the set $D$ is a FD-set of $G$, and so
$\fd(G) \le n - |X| + 1 \le n - c\log n + 1 = n - c^*\log n$ for some
constant~$c^*$.~\qed

\begin{prop}
If $G$ is a connected graph on $n \ge 6$ vertices satisfying $\fd(G)
= n - 2$, then $2 \le \rep(G) \le 4$. \label{prop5}
\end{prop}
\textbf{Proof.} A well-known elementary exercise states that every
graph of order at least two has two vertices with the same degree,
and so $\rep(G) \ge 2$. Hence it suffices for us to prove that
$\rep(G) \le 4$. Assume for the sake of contradiction, that $\rep(G)
\ge 5$. By Ramsey Theory, $r(K_3,K_3) = 6$. Further if $G$ is a graph
on five vertices such that neither $G$ nor its complement $\barG$
contains a copy of $K_3$, then $G = C_5$. Thus since $\rep(G) \ge 5$
there are either three vertices of the same degree in $G$ that induce
an independent set or clique in $G$ or there are five vertices of the
same degree in $G$ that induce a $C_5$. Let $X$ be the vertex set of
such an independent set or clique in $G$ of cardinality~$3$, if it
exists; otherwise let $X$ be the vertex set of such an induced
$5$-cycle in $G$. Then, $G[X]$ is $s$-regular with $s = 2$ or $s =
0$. Since $G$ is connected and the vertices in $X$ have all the same
degree in $G$, they also must have the same number of neighbors in $V
\setminus X$ (which is non-empty since $n \ge 6$), implying that
$\fd(G) \le n - 3$ or $\fd(G) \le n - 5$ in case case $G[X] =
C_5$. In all cases, we have $\fd(G) < n - 2$, a contradiction.
Therefore, $\rep(G) \le 4$.~\qed




\medskip
We remark that the restriction on the order $n \ge 6$ in the
statement of Proposition~\ref{prop5} is necessary since $G = C_5$ has
order~$n = 5$ and satisfies $\rep(G) = 5$ and $\fd(G) = 3 = n-2$.
Both Proposition~\ref{prop4} and Proposition~\ref{prop5} give further
evidence that a connected graph of large order achieving the upper
bound in Theorem~\ref{thm_n2} is highly non-regular. Further,
Proposition~\ref{prop4} gives a better bound than
Proposition~\ref{prop2}(a) when $r \ge \sqrt{n/\log n}$ and a better
bound than Proposition~\ref{prop2}(b) when $\chi(G) \ge n/\log n$.

\subsection{Trees}

In this section, we focus our attention on trees. We shall need the
following notation. A vertex of degree one is called a \emph{leaf}
and its neighbor is called a \emph{support vertex}. The set of
support vertices of a tree $T$ is denoted by $S_T$, while the set of
leaves by $L_T$. A neighbor of a vertex $v$ that is a leaf we call a
\emph{leaf-neighbor} of $v$. A \emph{strong support vertex} is a
vertex adjacent to at least two leaves. The \emph{corona} of a graph
$H$, denoted by $\coro(H)$, is the graph of order $2|V(H)|$ obtained
from $H$ by attaching a leaf to each vertex of~$H$. We note that
every vertex of $\coro(H)$ is a leaf or is a support vertex with
exactly one leaf-neighbor.

If we restrict our attention to trees, then the bound in
Theorem~\ref{thm_n2} can be improved significantly. For this purpose,
we recall that a  classical result of Ore~\cite{o9} established that
if $G$ is a graph of order $n$ with no isolated vertex, then
$\gamma(G) \le n/2$.
Further, Payan and Xuong~\cite{px9} showed that the only connected
graphs achieving equality in this bound are the $4$-cycle $C_4$ and
the corona $\coro(H)$ for a connected graph~$H$. We next establish an
upper bound on the fair domination number of a tree and characterize
the extremal trees. We begin with the following two observations.

\begin{ob}
Every 1FD-set in a graph contains all its strong support vertices.
\label{o:tree1}
\end{ob}
\proof Let $G$ be a graph and let $D$ be a 1FD-set in $G$. Let $v$ be
an arbitrary strong support vertex in $G$. If $v \notin D$, then in
order to dominate the leaf-neighbors of $v$, every leaf-neighbor of
$v$ belongs to $D$. Since $v$ has at least two leaf-neighbors, this
implies that $|N(v) \cap D| \ge 2$, a contradiction. Hence, $v \in
D$.~\qed

\begin{ob}
If $T$ is the corona of a tree and $T$ has order~$n$, then $\fd(T) =
n/2$. Further, $V(T)$ can be partitioned into two $\fd(T)$-sets.
\label{o:tree2}
\end{ob}
\proof The result is trivial for $n = 2$. Hence we may assume that
$T$ is the corona of a tree and $n \ge 4$. Then, $\gamma(T) = n/2$
 and we note that $|S_T| = |L_T| = n/2$. Further both
sets $S_T$ and $L_T$ form a 1FD-set of $T$, and so $\fd(T) \le
\ofd(T) \le n/2$. Since every FD-set of $T$ is a dominating set, we
have that $n/2 = \gamma(T) \le \fd(T) \le n/2$. Consequently, we must
have equality throughout this inequality chain. In particular,
$\fd(T) = n/2$ and both $S_T$ and $L_T$ are $\fd(T)$-sets.~\qed

\medskip
We are now in a position to prove the following result. In the proof,
we will deal with what we call \emph{$3$-end-paths}, which are paths
$xyz$ in a tree $T$ such that $x$ is a leaf, $N(y) = \{x,z\}$ and
$d_T(z) \ge 2$. We will call $z$ the \emph{base vertex} of the
$3$-end-path. Note that, since every tree has at least two leaves,
the corona of a tree on at least three vertices has at least two
$3$-end-paths sharing at most their base vertices.

\begin{thm}\label{n/2}
If $T$ is a tree of order~$n \ge 2$, then $\ofd(T) \le n/2$ with
equality if and only if $T$ is the corona of a tree. \label{thm1}
\end{thm}
\textbf{Proof.} By Observation~\ref{o:tree2} if $T$ is the corona of
a tree, then $\fd(T) = n/2$ and both $S_T$ and $L_T$ are 1FD-sets. We
will prove the statement by induction on $n$. If $2 \le n \le 7$,
this follows directly by checking all possible trees. This
establishes the case cases. For the inductive hypothesis, let $n \ge
8$ and assume that every tree $T'$ of order~$n'$, where $2 \le n' <
n$, satisfies $\ofd(T) \le n'/2$, with equality only if $T'$ is the
corona of a tree. Let $T$ be a tree of order~$n$.
If $T$ is a star $K_{1,n}$, then the central vertex of $T$ is a
1FD-set, implying that $\fd(T) = 1 < n/2$ and we are done. Hence we
may assume that $T$ is not a star. Therefore, $T$ contains a vertex
$w$ all of whose neighbors except for one, say $y$, are leaves. Let
$t$ be the number of leaf-neighbors of $w$, and so $t = d(w) - 1$. We
distinguish the following cases.

\medskip
\indent \emph{Case~1. Suppose that $t = 1$}. Then, $d(w) = 2$. Let
$z$ be the leaf-neighbor of $w$. Then, $ywz$ is a $3$-end-path in $T$
with $y$ as its base vertex. Let $T^* = T - \{w,z\}$ have
order~$n^*$, and so $n^* = n - 2$. Applying the inductive hypothesis
to $T^*$, $\ofd(T^*) \le n^*/2 = n/2 - 1$, with equality only if
$T^*$ is the corona of a tree. Let $D^*$ be a $\ofd(T^*)$-set.
Suppose first that $T^*$ is not the corona of a tree. Then, $|D^*| <
n/2 - 1$. Moreover, $T$ is not the corona of a tree. If $y \in D^*$,
then let $D = D^* \cup \{w\}$. If $y \notin D^*$, then let $D = D^*
\cup \{z\}$. In both cases $D$ is a 1FD-set of $T$, and so $\ofd(T)
\le |D| < (n/2 -1) + 1 = n/2$, and we are done. Hence we may assume
that $T^*$ is the corona of a tree.
Then, $y$ is either a leaf or a support vertex of $T^*$.

Suppose that $y$ is a leaf of $T^*$. Since $T^*$ has at least five
vertices, it contains at least two $3$-end-paths that are vertex
disjoint or that have at most their base vertices in common.
Therefore, $T^*$ contains a $3$-end-path, say $abc$ where $c$ is the
base vertex of the path, that has no vertex from $N[y]$. We now
consider the tree $T^{**} = T - \{a,b\}$. In $T^{**}$ we note that
the vertex $y$ has degree~$2$ and has no leaf neighbor. Hence,
$T^{**}$ is not the corona of a tree. Applying the inductive
hypothesis to $T^{**}$, $\ofd(T^{**}) < n/2 - 1$. As above, every
$\ofd(T^{**})$ can be extended to a 1FD-set of $T$ by adding to it
either $w$ or $z$, implying that $\ofd(T) < n/2$ and we are done.
Hence we may assume that $y$ is a support vertex of $T^*$. Then, $T$
is also a corona of a tree and hence, by Observation \ref{o:tree2},
$\fd(T) = n/2$ and we are done.

\medskip
\indent \emph{Case~2. Suppose that $t \ge 2$}. Then evidently $T$ is
not the corona of a tree. Let $x_1, x_2, \ldots x_t$ be the $t$
neighbor leaves of $w$ and let $T^* = T - \{x_1, x_2, \ldots, x_t\}$.
Suppose first that $T^*$ is not the corona of a tree. Applying the
induction hypothesis to $T^*$, we have that $\fd(T^*) < (n - t)/2$.
Let $D^*$ be a $\ofd(T^*)$-set. If $y \in D^*$, then let $D = D^*
\cup \{w\}$. If $y \notin D^*$, then $w \in D^*$ and let $D = D^*$.
In both cases, $D$ is a 1FD-set of $T$, and so $\fd(T) \le |D| \le
|D^*| + 1 < (n - t)/2 + 1 \le n/2$, and we are done. Hence we may
assume that $T^*$ is the corona of a tree. By
Observation~\ref{o:tree2}, $\fd(T^*) = (n - t)/2$ and $L_{T^*}$ is a
$\ofd(T^*)$-set. Since $w \in L_{T^*}$, it follows that $L_{T^*}$ is
also a 1FD-set of $T$, implying that $\fd(T) \le (n-t)/2 < n/2$, and
we are done.~\qed

\medskip
As an immediate consequence of Observation~\ref{o:tree2} and
Theorem~\ref{thm1}, we have the following result.

\begin{cor}
If $T$ is a tree of order~$n \ge 2$, then $\fd(T) \le n/2$ with
equality if and only if $T$ is the corona of a tree. \label{c:thm1}
\end{cor}

Recall that the \emph{$k$-domination number} $\gamma_k(G)$ of a graph
$G$ is the cardinality of a minimum $k$-dominating set, i.e., a set
$D$ of vertices such that every vertex outside $D$ has at least $k$
neighbors in $D$. In~\cite{FiJa}, Fink and Jacobson show that
$\gamma_2(T) \ge \lceil (n+1)/2 \rceil$ holds for every tree $T$ on
$n$ vertices. This fact and Theorem~\ref{n/2} allow us to prove that
any minimum FD-set in a tree is a 1FD-set.

\begin{ob}\label{fd_tree}
In a tree, every minimum FD-set is a 1FD-set.
\end{ob}
\textbf{Proof.} Let $T$ be a tree and let $D$ be a $\fd(T)$-set.
Then, $D$ is a kFD-set for some $k \ge 1$. Suppose that $k \ge 2$.
Then, $|N(x) \cap D| = k \ge 2$ for all $x \in V \setminus D$,
implying that $D$ is a $2$-dominating set of $T$ and thus $|D| \ge
\gamma_2(T)$. Hence by the result of Fink and Jacobson, we have that
$\fd(T) = |D| \ge \gamma_2(T) \ge \lceil (n+1)/2 \rceil > n/2$,
contradicting Theorem~\ref{n/2}. Hence, $k = 1$, and so $D$ is a
1FD-set.~\qed

\medskip
The set of non-leaves in a tree is a FD-set in the tree, implying the
following observation.

\begin{ob}
If $T$ is a tree on $n \ge 3$ vertices with $\ell$ leaves, then
$\fd(T) \le n - \ell$. \label{o:tree}
\end{ob}

We remark that if a tree has more leaves than internal vertices, then
the upper bound on the fair domination number of a tree given by
Observation~\ref{o:tree} is better than the upper bound of
Theorem~\ref{n/2}. The next theorem characterizes the trees where the
set of non-leaves is not a minimum FD-set. If $H$ is a subtree of a
tree $T$ such that $H$ is the corona of a tree, $H \ne T$ and $N(x)
\cap S_H = \emptyset$ for every vertex $x \in V \setminus V(H)$, then
we call $H$ a \emph{special corona-subtree} of $T$.

\begin{thm}
Let $T = (V,E)$ be a tree on $n \ge 3$ vertices with $\ell$ leaves.
Then the following assertions are equivalent: \\
\indent {\rm (i)} $\fd(T) < n - \ell$. \\
\indent {\rm (ii)} For every $\fd(T)$-set $D$, there are two
     vertices in $V \setminus D$ that are adjacent. \\
\indent {\rm (iii)} The tree $T$ contains a special corona-subtree.
\label{char:tree}
\end{thm}
\proof We show that $(i) \Rightarrow (ii) \Rightarrow (iii)
\Rightarrow (i)$. Let $D$ be an arbitrary $\fd(T)$-set. By
Observation~\ref{fd_tree}, the set $D$ is a 1FD-set in $T$.

$(i) \Rightarrow (ii)$: Suppose that $\fd(T) < n - \ell$. Then, $D
\ne V \setminus L_T$ and thus either $D \subset V \setminus L_T$ (but
$D \ne V \setminus L_T$) or $L_T \cap D \ne \emptyset$. Suppose $D
\subset V \setminus L_T$ and let $x \in V \setminus (D \cup L_T)$.
Since $d_T(x) \ge 2$ and $|N(x) \cap D| = 1$, there is a vertex $y
\in N(x) \setminus D$ and (ii) holds. Hence we may assume that $L_T
\cap D \ne \emptyset$. Let $z \in L_T \cap D$ and let $N(z) = \{y\}$.
If $y \in D$, then $D \setminus \{z\}$ is a 1FD-set in $T$,
contradicting the minimality of $D$. Hence, $y \notin D$, implying
that $N(y) \cap D = \{z\}$. Since $n \ge 3$, there is a vertex $x \in
N(y) \setminus D$ and, once again, (ii) holds.

$(ii) \Rightarrow (iii)$: Suppose that for every $\fd(T)$-set $S$
there are two vertices in $V \setminus S$ that are adjacent. Let $x$
and $y$ be two adjacent vertices in $V \setminus D$. We note that $x,
y \in V \setminus L_T$. Let $x' \in N(x) \cap D$ and $y' \in N(y)
\cap D$, and consider the subtree $T'= T[\{x,x',y,y'\}]$ of $T$.
Then, $L_{T'} = V(T') \cap D = \{x,y\}$ and $S_{L'} = V(H) \setminus
L(T') = \{x,y\}$. Among all subtrees of $T$ containing $x$ and $y$
with the property that the leaves of the subtree are vertices of $D$,
let $H$ be one of maximum order. Thus, $\{x,y\} \subset V(H)$, $L_H =
V(H) \cap D$ and $H$ is such a tree of maximum order. Let $S = V(H)
\setminus L_H$.

We show that $H$ is the corona of a tree. It suffices to show that
$|N(u) \cap L_H| = 1$ for all $u \in S$. Let $u \in S$ and suppose
first that $N(u) \cap L_H = \emptyset$. Since $L_H = V(H) \cap D$,
there is a vertex $y \in V \setminus V(H)$ such that $y \in D$ and $y
\in N(u)$. But then adding the vertex $y$ and the edge $uy$ to the
tree $H$ produces a tree $H'$ such that $\{x,y\} \subset V(H')$ and
$L_{H'} = V(H') \cap D$, contradicting the maximality of $H$. Hence,
$|N(u) \cap L_H| \ge 1$ for all $u \in S$. If $|N(u) \cap L_H| > 1$
for some $u \in S$, then $|N(u) \cap D| > 1$, contradicting the fact
that $D$ is a 1FD-set. Hence, $|N(u) \cap L_H| = 1$ for all $u \in
S$, implying that $H$ is the corona of a tree and $S_H=S$. If $T$ is the corona
of a tree, then by Observation~\ref{o:tree2}, the set $S_T$ is a
$\fd(T)$-set that does not satisfy (ii), a contradiction. Hence, $T$
is not the corona of a tree. Thus, $H \ne T$.

Let $z \in V \setminus V(H)$ and suppose that there is a vertex $v
\in N(z) \cap S$. Let $w$ be the vertex of $D$ that is adjacent to
$z$ in $T$. Since $T$ is cycle-free, $w \in V \setminus V(H)$. But
then adding the vertex $z$ and the edge $wz$ to the tree $H$ produces
a tree $H'$ such that $\{x,y\} \subset V(H')$ and $L_{H'} = V(H')
\cap D$, contradicting the maximality of $H$. Hence, $N(z) \cap S =
\emptyset$. This is true for all vertices $z \in V \setminus V(H)$.
Consequently, $H$ is a special corona-subtree of $T$, and so (iii)
holds.

$(iii) \Rightarrow (i)$: Suppose that the tree $T$ contains a special
corona-subtree $H$. Let $D = (V \setminus (L_T \cup S_H)) \cup (L_H
\cap L_T)$. Let $v \in V \setminus D = S_H \cup (L_T \setminus
V(H))$. If $v \in S_H$, then since $H$ is a special corona-subtree of
$T$ and $L_H \subseteq D$, we have that $|N_T(v) \cap D| = |N_T(v)
\cap L_H| = 1$. If $v \in L_T \setminus V(H)$, then since the
neighbor of $v$ in $T$ belongs to $D$, we once again have that
$|N_T(v) \cap D| = 1$. Hence, $D$ is a 1FD-set of $T$, and so $\fd(T)
\le |D| = |V \setminus L_T| - |S_H| + |L_H \cap L_T|$. Since $H$ is a
special corona-subtree of $T$,  at least one leaf of $T$ is not a
leaf of $H$. Therefore, $|L_H \cap L_T| < |L_H| = |S_H|$, implying
that $\fd(T) < |V \setminus L_T| = n - \ell$, and so (i) holds.~\qed

\medskip
As an immediate consequence of Theorem\ref{char:tree}, we obtain the
following characterization of the trees whose set of non-leaves is a
minimum FD-set.

\begin{cor}
Let $T = (V,E)$ be a tree on $n \ge 3$ vertices with $\ell$ leaves.
Then the following assertions are equivalent: \\
\indent {\rm (i)} $\fd(T) = n - \ell$. \\
\indent {\rm (ii)} There is a $\fd(T)$-set $D$ such that
     $V \setminus D$ is an independent set. \\
\indent {\rm (iii)} The tree $T$ contains no special corona-subtree.
\label{char:treec}
\end{cor}

\subsection{Maximal Outerplanar Graphs}

A maximal outerplanar graph, abbreviated \emph{MOP}, is a
triangulation of the polygon. It is well-known that every bounded
face of a MOP is a triangle. Further a MOP on $n$ vertices is
$3$-colorable, $2$-degenerate (i.e., every subgraph contains a vertex
with degree at most two), has exactly $2n-3$ edges and the
neighborhood of every vertex of the graph induces a path (see
\cite{HoSt, KoWe}). In particular it follows that every MOP $G$ on
$n$ vertices has $\delta(G) = 2$ and that the vertices of degree $2$
are independent when $n \ge 4$. Note also that $\bard(G) = 4 - 6/n$
and hence from Proposition \ref{prop2}(b) it follows that $\fd(G) \le
14n/15$. With more effort, this bound can be improved to $17n/19$, as
we will show in the next theorem. But first, we need to prove the
following lemma.

\begin{lem}\label{bip}
If $G$ is a MOP on $n \ge 3$ vertices, then the vertices of
degree~$3$ induce a bipartite graph.
\end{lem}
\proof We proceed by induction on $n$. If $n = 3$ or $n=4$, then the
theorem holds trivially. This establishes the base cases. Let $n \ge
5$ and assume that in every MOP of order~$n'$, where $3 \le n' < n$,
the vertices of degree~$3$ induce a bipartite graph. Let $G$ be a MOP
on $n$ vertices. Since $G$ is maximal outerplanar, there is a vertex
$x$ of degree~$2$ whose deletion results in a graph $G^*$ on $n-1$
vertices which is again a MOP. Let $u$ and $v$ be the neighbors of
$x$ in $G$. Then, $u$ and $v$ are adjacent. Further since the
vertices of degree~$2$ in $G$ are independent, we have that $d_G(u)
\ge 3$ and $d_G(v) \ge 3$. If $d_G(u) = 3$ and $d_G(v) = 3$, then $u$
and $v$ would be two adjacent vertices of degree~$2$ in the MOP $G^*$
of order at least~$4$, which is not allowed. Hence, renaming $u$ and
$v$ if necessary, we may assume that $d_G(u) \ge 4$ and $d_G(v) \ge
3$.
Let $B$ and $B^*$ be the sets of vertices of degree $3$ in $G$ and
$G^*$, respectively. By the induction hypothesis, the set $B^*$
induces a bipartite graph. If $d_G(v) \ge 4$, then $B = B^* \setminus
\{u,v\}$ and $B$ induces also a bipartite graph in $G$, as desired.
Hence we may assume that $d_G(v) = 3$. If $d_G(u) \ge 5$, then $u, v
\notin B^*$ and let $B = B^* \cup \{v\}$. If $d_G(u) = 4$, then let
$B = (B^* \setminus \{u\}) \cup \{v\}$. In both cases, the vertex $v$
is adjacent to at most one vertex of degree~$3$ in $G$, say $z$, and
so $v$ can be added to a partite set of $G^*[B^*]$ that does not
contain $z$, showing then that $G[B]$ is bipartite.~\qed

\medskip
We are now in a position to present the following upper bound on the
fair domination number of a MOP.

\begin{thm}
If $G$ is a MOP on $n \ge 3$ vertices, then $\fd(G) < 17n/19$.
\end{thm}
\proof For $n = 3$, then $\fd(G) = 1 < 17n/19$. Hence we may assume
that $n \ge 4$. For $i = 2,3,4,5$, let $V_i$ denote the set of
vertices of degree~$i$ in $G$ and let $|V_i| = n_i$. Moreover let $t$ be
the number of vertices of degree at least~$6$ in $G$. Then, $n_2 +
n_3 + n_4 + n_5 + t = n$ and $2 n_2 + 3 n_3 + 4 n_4 + 5 n_5 + 6 t \le
2 m = 2 (2n - 3) = 4n - 6$, counting first vertices and then edges.
If $n_2 \le 2n/19$, $n_3 \le 4n/19$, $n_4 \le 6n/19$ and $n_5 \le
6n/19$, then
\begin{eqnarray*}
4n-6 &\ge& 2 n_2 + 3 n_3 + 4 n_4 + 5 n_5 + 6 t \\
         &=& 2 n_2 + 3 n_3 + 4 n_4 + 5 n_5 + 6  (n - ( n_2 + n_3 + n_4 +  n_5 ))\\
         &=& 6n - 4 n_2 - 3 n_3 - 2 n_4 - n_5 \\
        &\ge& 6n - \frac{8}{19} n - \frac{12}{19} n - \frac{12}{19} n - \frac{6}{19} n \\
        & =&  4n,
\end{eqnarray*}
a contradiction. Hence, $n_2 > 2n/19$, $n_3 > 4n/19$, $n_4 > 6n/19$
or $n_5 > 6n/19$. Suppose first that $n_2 > 2n/19$. Then, as $n \ge
4$, the set $V_2$ is independent and is therefore an out-regular set,
and so $\outr(G) \ge n_2 > 2n/19$. Hence by Proposition~\ref{prop1},
we have that $\fd(G) = n - \outr(G) < 17n/19$. Suppose $n_3 > 4n/19$.
Let $V_3'$ be a maximum independent subset of $V_3$. By
Lemma~\ref{bip}, the graph $G[V_3]$ is a bipartite graph, and so
$|V_3'| \ge |V_3|/2 > 2n/19$. Since $V_3'$ is an out-regular set, we
have that $\outr(G) \ge |V_3'| > 2n/19$. Hence by
Proposition~\ref{prop1}, we have that $\fd(G) < 17n/19$. Suppose $n_4
> 6n/19$. Let $V_4'$ be a maximum independent subset of $V_4$. Since
$G$ is $3$-colorable, we have that $|V_4'| \ge |V_4|/3 > 2n/19$. The
set $V_4'$ is an out-regular set, and so $\outr(G) \ge |V_4'| >
2n/19$, implying by Proposition~\ref{prop1} that $\fd(G) < 17n/19$.
Analogously if $n_5 > 6n/19$, then $\fd(G) < 17n/19$. In all four
cases, we have $\fd(G) < 17n/19$, as desired.~\qed

\subsection{Nordhaus-Gaddum-Type Bounds}

In this section we consider Nordhaus-Gaddum-type bounds for the fair
domination number of a graph.

\begin{thm} \label{t:NordG}
Let $G$ be a graph on $n$ vertices. Then the following holds.\\
\indent {\rm (a)} If $n \ge 5$, then $3 \le \fd(G) + \fd(\barG) \le
2n-4$ and both bounds are sharp.\\
\indent {\rm (b)} If $n \ge 4$, then $2 \le \fd(G) \cdot \fd(\barG)
\le (n-2)^2$ and both bounds are sharp.
\end{thm}
\proof We will first prove both upper bounds. Without loss of
generality, we may assume that $G$ is connected, and so, by
Theorem~\ref{thm_n2}, $\fd(G) \le n - 2$. If $\barG$ has size
$m(\barG) \ge 2$, then by Corollary~\ref{c:bound}, $\fd(\barG) \le n
- 2$, and so $\fd(G) + \fd(\barG) \le 2n-4$ and $\fd(G) \cdot
\fd(\barG) \le (n-2)^2$, as desired. If $m(\barG) = 1$, then $\barG
= \overline{K}_{n-2} \cup K_2$ and $G = K_n - e$ for some edge
$e$ of $K_n$. Hence, $\fd(G) = 1$ and $\fd(\barG) = n-1$, implying
$\fd(G) + \fd(\barG) = n \le 2n -4$ and $\fd(G) \cdot \fd(\barG) = n
-1 \le (n-2)^2$ when $n \ge 4$. Finally, if $m(\barG) = 0$, then
$\barG = \overline{K_n}$ and $G = K_n$. In this case we have
${\rm fd}(G) + \fd(\barG) = 1 + n \le 2n-4$, for $n \ge 5$, and ${\rm
fd}(G) \cdot \fd(\barG)  = n \le (n-2)^2$, for $n \ge 4$. In all
cases, $\fd(G) + \fd(\barG) \le 2n-4$, for $n \ge 5$, and $\fd(G)\,
\fd(\barG) \le (n-2)^2$, for $n \ge 4$.

That the upper bounds are sharp may be seen as follows. For $n \ge 6$
even, let $G = H_{n/2}$, while for $n \ge 7$ odd, let $G =
F_{(n-1)/2}$ where the graphs $H_n$ and $F_n$ are defined as in the
proof of Theorem~\ref{thm_n2}. Then, $G$ has order~$n$ and, as shown
in the proof of Theorem~\ref{thm_n2}, $\fd(G) = n - 2$. Since the
graph $\barG$ is connected, by Theorem~\ref{thm:Gbar}(a), we have
that $\fd(\barG) = n - 2$, and so $\fd(G) + \fd(\barG) = 2n-4$ and
$\fd(G)\, \fd(\barG) = (n-2)^2$. For $n = 4, 5$, we simply take $G =
C_n$, which satisfies $\fd(G) = \fd(\barG) = n - 2$.

For the lower bounds, we may assume, without loss of generality, that
$\fd(G) \le \fd(\barG)$. If $\fd(G) \ge 2$, then $\fd(G) + \fd(\barG)
\ge 3$ and $\fd(G) \cdot \fd(\barG) \ge 2$. Hence we may assume that
$\fd(G) = 1$. But then $G$ has a vertex of degree~$n-1$, implying
that $\barG$ is not connected. Hence, $\fd(\barG) \ge 2$ and we
obtain $\fd(G) + \fd(\barG) \ge 3$ and $\fd(G) \cdot \fd(\barG) \ge
2$. This establishes the desired lower bounds. That the lower bounds
are sharp may be seen by taking $G = K_{1,n-1}$.~\qed

\medskip
We remark that if $G = K_4$, then $G$ has order~$n = 4$ and $\fd(G) +
\fd(\barG) = 1 + 4 = 5 > 2n-4$. Hence the constraint on the order~$n
\ge 5$ in the statement of Theorem~\ref{t:NordG}(a) cannot be
relaxed. Moreover neither can the condition $n \ge 4$ in
Theorem~\ref{t:NordG}(b) be relaxed since if $G = K_3$, then $G$ has
order $n=3$ and we have $\fd(G) = 1$ and $\fd(\barG) = 3$, implying
that $\fd(G)\, \fd(\barG) = 3 > (n-2)^2 = 1$.

\subsection{Unions of Graphs}

In this section we investigate the fair domination of disjoint unions
of graphs.
We shall prove.

\begin{thm}
Let $G_1,\ldots,G_k$ be $k \ge 1$ graphs. Let $H$ be the disjoint
union $\bigcup_{i=1}^k G_i$ of $G_1,\ldots,G_k$ and let $H$ have
order~$n$. Then the following holds. Then,
\[
\fd(H) - \sum_{i=1}^k \fd(G_i) \le \frac{1}{k} (k-1)(n-k),
\]
and this bound is sharp for $k = 1,2$.
 \label{union}
\end{thm}
\textbf{Proof.} For $k = 1$ the result is trivial since both sides of
the inequality are zero. Hence we may assume that $k \ge 2$. Let $n_i = |V(G_i)|$ for $1 \le i \le k$. Renaming the graphs $G_1,\ldots,G_k$, if necessary, we may assume that $n_1 \le n_2 \le \cdots \le n_k$, and so $n_k \ge n/k$ and $n -
n_k \le (k-1)n/k$. Let $D_k$ be a $\fd(G_k)$-set. Then, $D_k \cup (V(H) \setminus V(G_k))$ is a FD-set of $H$ and thus
\[
\begin{array}{lcl}
\fd(H) - \sum_{i=1}^k \fd(G_i) & \le & \displaystyle{ |D_k| + n - n_k - \sum_{i=1}^k \fd(G_i)  } \\ \2
& = &  \displaystyle{ n - n_k - \sum_{i=1}^{k-1} \fd(G_i)    }  \\   \2
& \le &  \displaystyle{ n - n_k - (k-1)}\\ \2
& \le & \displaystyle{ \frac{1}{k} (k-1)(n-k).}
\end{array}
\]

This establishes the desired upper bound.

We show next that the bound is sharp for $k = 2$. For $n \ge 3$, let
$G_1 = K_{1,2n-1}$ and let $G_2 = K_{2,n-1,n-1}$. Let $H = G_1 \cup
G_2$, and so $H$ has order~$4n$.

We note that there are only two values of $\ell$ for which $G_1$
has an $\ell$FD set $D$ different from $V(G_1)$. The first value of $\ell$ is $\ell =
1$ when $D$ consists of the central vertex of $G_1$ and any set of
$t$ leaves, where $0 \le t \le 2(n-1)$. The second value of $\ell$ is
$\ell = 2n-1$ when $D$ consists of all leaves in $G_1$. Hence, $\ell \in \{1,2n-1\}$. In particular, we
note that $\fd(G_1) = 1$.

Next we consider the graph $G_2$. Since $G_2$ has no dominating
vertex adjacent to every other vertex, we note that $G_2$ has no
$\ell$FD-set for $\ell \in \{1,2n-1\}$. In particular, $\fd(G_2) \ge
2$. However the partite set of cardinality~$2$ is a FD-set in $G_2$,
and so $\fd(G_2) \le 2$. Consequently, $\fd(G_2) = 2$.

From our earlier observations, a FD-set of $H$ is formed by either
taking a FD-set in $G_1$ and all vertices in $G_2$ or by taking a
FD-set in $G_2$ and all vertices in $G_1$. Therefore, $\fd(H) \le
\min \{\fd(G_1) + |V(G_2)|, \fd(G_2) + |V(G_1)|\} = 2n+1$. Thus,
\[
\fd(H) - \sum_{i=1}^k \fd(G_i) = 2(n-1) = \frac{1}{2} (4n-2) = \frac{1}{k} (k-1)(|V(H)|-k),
\]
realizing the upper bound for $k = 2$.~\qed

\medskip
We remark that the fair domination number is highly sensitive with
respect to edge deletion or edge addition. For $n \ge 3$, let $G_1 =
K_{1,2n-1}$ and let $G_2 = K_{2,n-1,n-1}$. Let $v_1$ denote the
central vertex of $G_1$ and let $\{u_2,v_2\}$ be the partite set of
$G_2$ of cardinality~$2$. Let $H = G_1 \cup G_2$, and so $H$ has
order~$4n$, and let $G = H + u_2v_2$. As shown in the proof of
Theorem~\ref{union}, we have that $\fd(H) = 2n+1$. However the set
$\{v_1,v_2\}$ is a FD-set of $G$ since every vertex in $V(G)
\setminus \{v_1,v_2\}$ is adjacent to exactly one of $v_1$ and $v_2$.
Hence, $\fd(G) = 2$. Therefore we have the following observation.

\begin{prop}
There exists infinitely many graph $G$ such that $\fd(G + e) - \fd(G)
\ge |V(G)|/2 - 1$ for some edge $e \in E(\barG)$.
 \label{sensitive}
\end{prop}

\section{Closing Remarks and Open Questions}

This paper, in which we introduce the notion of fair domination,
suggests many possible direction of further research. We close with
the following list of open problems that we have yet to settle.

\begin{prob}
Find a polynomial time algorithm to compute $\fd(T)$ for trees $T$.
\end{prob}

\begin{prob}
Improve upon the bounds of Proposition~\ref{prop2}. In particular,
find best possible upper bounds for MOP's, maximum planar graphs and
regular graphs.
\end{prob}

\begin{prob}
Find $\fd(G)$ for other families of graphs $G$, in particular the
grid $P_n \times P_m$ and the torus $C_n \times C_m$.
\end{prob}

\begin{prob} Is it true that Theorem~\ref{union} is sharp also for $k \ge 3$?
 \end{prob}

\begin{prob} In view of Proposition \ref{sensitive}, find the exact value of $\max |\fd(G+ e) - \fd(G)|$ over all graphs $G$ on $n$ vertices.
 \end{prob}

\begin{prob} Characterize the graphs $G$ on $n$ vertices for which $\fd(G) = n-2$.
 \end{prob}

\begin{prob} What can be said about the fair domination number $\fd(G)$ of the random graph $G(n,p)$?
\end{prob}


\medskip
\begin{thebibliography}{99}

\bibitem{AlKrSu} N. Alon, M. Krivelevich, and B. Sudakov, Large
    nearly
    regular induced subgraphs. \textit{SIAM J. Discrete Math.}
    \textbf{22} (2008), no. 4, 1325--1337.

\bibitem{CKV} D. M. Cardoso and M. Kami\`{n}ski, V. Lozin, Maximum
    $k$-regular induced subgraphs. \textit{J. Comb. Optim.}
    \textbf{14} (2007), no. 4, 455--463.

\bibitem{Ca79} Y. Caro, New results on the independence number.
    \textit{Tech. Report}, \emph{Tel-Aviv University} (1979).

\bibitem{CaWe09}  Y. Caro and D. B. West, Repetition number of
    graphs. \textit{Electronic J.  Combin.} \textbf{16} (2009) \#R7.

\bibitem{Cs07} B. Csaba, Regular spanning subgraphs of bipartite
    graphs of high minimum degree. \textit{Electron. J. Combin.} {\bf
    14} (2007), no. 1, Note 21, 7 pp. (electronic).

\bibitem{HoVo} A. Hoffmann and L. Volkmann, On regular factors in
    regular graphs with small radius. \textit{Electron. J. Combin.}
    \textbf{11} (2004), no. 1, Research Paper 7, 7 pp. (electronic).

\bibitem{FiJa} J. F. Fink and M. S. Jacobson, $n$-Domination in
    graphs.
    \textit{Graph Theory with Applications to Algorithms and Computer
    Science.} John Wiley and Sons. New York (1985), 282-300.

\bibitem{hhs1} T. W. Haynes, S. T. Hedetniemi, and P. J. Slater
    (eds), \emph{Fundamentals of Domination in Graphs}, Marcel
    Dekker, Inc. New York, 1998.

\bibitem{hhs2} T. W. Haynes, S. T. Hedetniemi, and P. J. Slater
    (eds), \emph{Domination in Graphs: Advanced Topics}, Marcel
    Dekker, Inc. New York, 1998.

\bibitem{HoSt} G. Hopkins and W. Staton, Outerplanarity without
    topology, Bull. Inst. Combin. Appl. \textbf{21} (1997), 112--116.

\bibitem{KoWe} A. Kostochka and D. B. West, Every
    outerplanar
    graph is the union of two interval graphs.  \textit{Proceedings
    of the Thirtieth Southeastern International Conference on
    Combinatorics, Graph Theory, and Computing} (Boca Raton, FL,
    1999).  Congr. Numer. \textbf{139} (1999), 5--8.

\bibitem{o9} O. Ore, \emph{Theory of graphs}. \textit{Amer. Math.
    Soc. Transl.} \textbf{38} (Amer. Math. Soc., Providence, RI,
    1962), 206--212.


\bibitem{px9} C. Payan and   N. H. Xuong, Domination-balanced graphs.
    \textit{J. Graph Theory} \textbf{6} (1982), 23--32.

\bibitem{Wei81} V. K. Wei, A lower bound on the stability number of a
    simple graph. \emph{Bell Lab. Tech. Memo.} No. 81-11217-9 (1981).



\end{thebibliography}
\end{document}